\documentclass[12pt,reqno,a4paper]{amsart}

\usepackage{amsmath,
amsthm,
cite}

\usepackage{txfonts}

\newtheorem{theorem}{Theorem}
\newtheorem{prop}[theorem]{Proposition}

\usepackage{titlesec}
\titleformat{\section}
{\bfseries\scshape\centering}
{\thesection.}{.5em}{}
\titleformat{\subsection}
{\rmfamily\bfseries}
{\thesubsection}{.5em}{}

\begin{document}

\renewcommand{\thefootnote}{\fnsymbol{footnote}}





\title{
  Difference Equation of the Colored Jones Polynomial for
  Torus Knot
}


    \author{Kazuhiro \textsc{Hikami}}


  \address{Department of Physics, Graduate School of Science,
    University of Tokyo,
    Hongo 7--3--1, Bunkyo, Tokyo 113--0033, Japan.
    }

    \email{hikami@phys.s.u-tokyo.ac.jp}
    \urladdr{http://gogh.phys.s.u-tokyo.ac.jp/{\textasciitilde}hikami/}

\vspace{18pt}
\date{March 11, 2004}

\begin{abstract}
We prove
that the $N$-colored Jones polynomial for the torus knot
$\mathcal{T}_{s,t}$ satisfies the second order difference equation,
which  reduces to the first order
difference equation for a case of $\mathcal{T}_{2,2m+1}$.
We show that the A-polynomial of the torus knot can be derived
from the
difference equation.
Also constructed is 
a $q$-hypergeometric type expression of  the colored
Jones polynomial for $\mathcal{T}_{2,2m+1}$.

\end{abstract}



\subjclass[2000]{
}


\maketitle
\section{Introduction}

The $N$-colored Jones polynomial
$J_\mathcal{K}(N)$
is one of the quantum invariants for
knot $\mathcal{K}$.
It is associated with the $N$-dimensional irreducible representation of
$sl(2)$,
 and is powerful to classify knots.
Motivated by ``volume conjecture''~\cite{Kasha96b,MuraMura99a}
saying that
a hyperbolic volume of the knot complement dominates an asymptotic behavior of 
the colored Jones polynomial $J_\mathcal{K}(N)$ at $q=\mathrm{e}^{2 \pi \mathrm{i}/N}$,
it
receives much interests toward  a geometrical and topological
interpretation of the quantum invariants.

Recently  another intriguing structure of the colored Jones polynomial was put forward;
it was shown that the $N$-colored Jones polynomial
$J_{\mathcal{K}}(N)$  can be
written in a $q$-hypergeometric form, and that it satisfies a
recursion relation with respect to $N$~\cite{GaroTQLe03a}.
It was further demonstrated 
for trefoil and figure-eight knot~\cite{SGarou03b} 
that  a recursion relation is related to the
A-polynomial $A_\mathcal{K}(L,M)$
(see also Ref.~\citen{FrohGelcLofa01a}),
which
denotes  an algebraic curve of eigenvalues of the $SL(2,\mathbb{C})$
representation of the boundary torus of knot $\mathcal{K}$~\cite{CCGLS94a}.
As the A-polynomial contains many geometrical informations such as the boundary
slopes of the knot, this ``AJ conjecture'' may help our geometrical
understanding of the colored Jones polynomial.

In this article, we study torus knot $\mathcal{T}_{s,t}$ where $s$ and
$t$ are relatively prime integers.
We  prove that the $N$-colored Jones polynomial
$J_\mathcal{K}(N)$
for the torus knots
$\mathcal{K}=\mathcal{T}_{s,t}$ satisfies the second order recursion
relation~\eqref{difference_Jones}
[Theorem~\ref{thm:difference_Jones}],
which reduces to the first order~\eqref{torus_difference}
only  in a case of
$\mathcal{K}=\mathcal{T}_{2,2m+1}$
[Theorem~\ref{thm:torus_difference}].
Furthermore we shall show that this difference operator gives the
A-polynomial of the torus knot as was demonstrated in
Ref.~\citen{SGarou03b}
[Theorem~\ref{thm:AJ_torus}].
Throughout this article, we normalize the colored Jones polynomial to be
\begin{equation*}
  J_{\text{unknot}}(N)=1
\end{equation*}


\section{\mathversion{bold}
Colored Jones Polynomial}

The Alexander polynomial
$\Delta_{\mathcal{K}}(A)$ for the torus knot
$\mathcal{K}=\mathcal{T}_{s,t}$ is known to be
(see \emph{e.g.}  Ref.~\citen{Licko97Book})
\begin{equation}
  \Delta_{\mathcal{T}_{s,t}}(A)
  =
  \frac{
    (A^{{1}/{2}} - A^{-{1}/{2}}) 
    (A^{s t/2} - A^{- s t / 2}) 
  }{
    (A^{{s}/{2} } - A^{-{s}/{2} }  ) 
    (A^{{t}/{2} } - A^{-{t}/{2}  }) 
  }
\end{equation}
We see that an inverse of the Alexander polynomial is expanded in $A\to\infty$
as
\begin{equation}
  \label{define_Chi}
  \frac{A^{{1}/{2}}-A^{-{1}/{2}}}{\Delta_{\mathcal{T}_{s,t}}(A)}
  =
  \sum_{n=0}^\infty \chi_{2 s t}(n) \, A^{-n/2}
\end{equation}
where $\chi_{2 s t}(n)$ is the periodic function with
modulus $2 \, s\, t$~\cite{KHikami03c};
\begin{equation*}
  \begin{array}{c|ccccc}
    n \mod 2 \, s \, t
    & s \, t - s - t & s \, t - t+s & s \, t + t -s & s \, t + s + t
    & \text{others}
    \\
    \hline
    \chi_{2 s t}(n) &
    1 & -1 & -1 & 1 & 0
  \end{array}
\end{equation*}
Using this periodic function,
we define the function $K_{{s,t}}(N)$ 
by
\begin{equation}
  \label{define_K}
  K_{s,t}(N) 
  =
  q^{\frac{1}{4} N \bigl( 2 ( s t -s -t) - s t N \bigr)}
  \sum_{k=0}^{s t N}
  \chi_{2 s t}(s \, t \, N -k) \,
  q^{\frac{k^2-(s t -s -t)^2}{4 s t}}
\end{equation}

\begin{theorem}
  \label{theorem:torus}
  The $N$-colored Jones polynomial for the torus knot $\mathcal{T}_{s,t}$
  is given by
  \begin{equation}
    \label{H_and_Jones}
    J_{\mathcal{T}_{s,t}}(N)
    =
    \frac{q^{\frac{1}{2} (s-1)(t-1) (1-N)}}{1-q^{-N}} \,
    K_{{s,t}}(N)
  \end{equation}
\end{theorem}

As we have defined the function $K_{{s,t}}(N)$ from an
expansion of the Alexander polynomial, this theorem reveals a
connection between the colored Jones polynomial and the Alexander
polynomial for a case of the torus knot.
It should be noted
that a relationship between these polynomials
is known  for arbitrary knot $\mathcal{K}$
based on  a slightly different expansion  as the
Melvin--Morton  conjecture~\cite{MelvMort95a,Rozan96c}, which was proved in
Ref.~\citen{BarNaGarou96a}.

To prove Theorem~\ref{theorem:torus}, we use a previously known result for the colored Jones polynomial.
\begin{prop}[\cite{Mort95a}]
  The $N$-colored Jones polynomial for the torus knot
  $\mathcal{T}_{s,t}$ is computed as
  \begin{equation}
    \label{Jones_Morton}
    J_{\mathcal{T}_{s,t}}(N)
    =
    \frac{q^{\frac{1}{4} s t (1-N^2)}}{q^{\frac{N}{2}} - q^{-\frac{N}{2}}}
    \sum_{r=-\frac{N-1}{2}}^{\frac{N-1}{2}}
    \left(
      q^{ s t r^2 - (s+t) r + \frac{1}{2}}
      -
      q^{ s t r^2 - (s-t) r - \frac{1}{2}}
    \right)
  \end{equation}
\end{prop}

\begin{proof}[Proof of Theorem~\ref{theorem:torus}]
  We  first prove for a case of  $N$ being even.
  We have
  \begin{align*}
    & \sum_{k=0}^{s t N}
    \chi_{2 s t}(s \, t \, N -k) \, q^{\frac{k^2- (s t - s -t )^2}{4 s t}}
    =
    \sum_{k=0}^{s t N}
    \chi_{2 s t}(k) \, q^{\frac{k^2- (s t - s -t )^2}{4 s t}}
    \\
    & =
    \sum_{k=0}^{\frac{N}{2} -1}
    \left(
      q^{k(s t k + s t -s-t)}
      -
      q^{(s k + s -1)( t k +1)}
      -
      q^{(s k + 1)(t k + t -1)}
      +
      q^{(k+1) ( s t k + s-1)}
    \right)
    \\
    & =
    \sum_{k=-\frac{N}{2}}^{\frac{N}{2} -1}
    \left(
      q^{k ( s t k + s t - s - t) }
      -
      q^{ (s k + s -1) (t k+1)}
    \right)
  \end{align*}
  In the second equality, we have
substituted a non-zero value of the periodic function $\chi_{2 s t}(k)$.
  In the last equality, we have replaced
  $k$ with $-k-1$ in both the third and the fourth terms in a parenthesis.
  By further replacing  $k$ with $k-\frac{1}{2}$ in the last expression
  and comparing with an
  explicit form given in eq.~\eqref{Jones_Morton},
  we get an assertion of the theorem.
  
  For a case of  $N$  being odd, we can prove it in a same manner.
\end{proof}


\begin{prop}
  Let the function $K_{s,t}(N)$ be defined by eq.~\eqref{define_K}.
  Then it satisfies the following difference equation;
  \begin{equation}
    \label{difference_K}
    K_{s,t}(N)
    =
    1 - q^{s(1-N)-1} - q^{t(1-N)-1} + q^{(s+t)(1-N)}
    +
    q^{s t(2-N) - s- t} \,
    K_{s,t}(N-2)
  \end{equation}
\end{prop}

\begin{proof}
  We decompose a sum in eq.~\eqref{define_K} into
  $\sum_{k=0}^{s t(N-2)} + \sum_{k=s t (N-2)}^{s t N}$.
  {}From the first sum  we obtain $K_{s,t}(N-2)$.
  The second sum can be written explicitly as what appeared  in
  eq.~\eqref{difference_K} using a property of $\chi_{2 s t}(k)$.
\end{proof}

As we have obtained  a relationship between the colored Jones polynomial
$J_{\mathcal{T}_{s,t}}(N)$ and $K_{s,t}(N)$,  it is straightforward to
obtain a following  theorem.

\begin{theorem}
  \label{thm:difference_Jones}
  The $N$-colored Jones polynomial for the torus knot
  $\mathcal{T}_{s,t}$ fulfills a recursion relation of the second
  order;
  \begin{multline}
    \label{difference_Jones}
    J_{\mathcal{T}_{s,t}}(N)
    =
    \frac{
      q^{\frac{1}{2} (s-1) (t-1) (1-N)}
    }{
      1-q^{-N}
    } \,
    \left(
      1 - q^{s (1-N)-1} - q^{t(1-N)-1} + q^{(s+t)(1-N)}
    \right)
    \\
    +
    \frac{1-q^{2-N}}{1-q^{-N}} \,
    q^{s t(1-N) - 1 } \,
    J_{\mathcal{T}_{s,t}}(N-2)
  \end{multline}
\end{theorem}

We should  note that this difference equation can be
directly  derived  by use of
Morton's expression~\eqref{Jones_Morton}.
A  benefit of  our expression~\eqref{H_and_Jones} is in reducing the
difference equation~\eqref{difference_Jones} into the first order
difference equation in a case of
$\mathcal{T}_{2,2m+1}$.

\begin{prop}[see \emph{e.g.} Ref.~\citen{KHikami02c}]
  The function $K_{{2,2m+1}}(N)$ satisfies the difference equation of the first order,
  \begin{equation}
    \label{difference_K_2}
    K_{2, 2m+1}(N)
    =
    1 - q^{1-2N}
    - q^{2m(1-N)-N} \,
    K_{2,2m+1}(N-1)
  \end{equation}
\end{prop}
\begin{proof}
  We note that
  the function $\chi_{8m+4}(n)$ has an anti-periodicity,
  $\chi_{8m+4}(n+4m+2) = - \chi_{8m+4}(n)$.
  Then a proof follows in a same method with that of
  eq.~\eqref{difference_K}.
\end{proof}

This proposition simplifies the difference equation of the colored Jones polynomial $J_{\mathcal{K}}(N)$
for $\mathcal{K}=\mathcal{T}_{2,2m+1}$ as follows;

\begin{theorem}
  \label{thm:torus_difference}
  The $N$-colored Jones polynomial for the torus knot $\mathcal{T}_{2,2m+1}$ solves the difference equation of the first order;
  \begin{equation}
    \label{torus_difference}
    J_{\mathcal{T}_{2,2m+1}}(N)
    = q^{m(1-N)} \, \frac{1-q^{1-2 N}}{1-q^{-N}}
    -
    q^{m-(2m+1)N} \, \frac{1-q^{1-N}}{1-q^{-N}} \,
    J_{\mathcal{T}_{2,2m+1}}(N-1)
  \end{equation}
\end{theorem}
This recursion relation  coincides with a result in Ref.~\citen{SGarou03b} for the trefoil
$m=1$
(we need to replace $q$ with $q^{-1}$).
We remark that
in
Ref.~\citen{GelcSain03a}
by  a different approach proposed was the difference equation of the
colored Jones polynomial for the
torus knot $\mathcal{T}_{2,2m+1}$, which is much involved.

\section{A-Polynomial}

In Ref.~\citen{SGarou03b} 
the ``AJ conjecture'' is proposed;
the  \emph{homogeneous} difference equation of the colored Jones polynomial
$J_\mathcal{K}(N)$ for knot $\mathcal{K}$
gives the A-polynomial $A_\mathcal{K}(L,M)$
for knot $\mathcal{K}$.
More precisely    we rewrite the difference equation of the
colored Jones polynomial,
$\sum_{k\geq 0} a_k \, J_{\mathcal{K}}(N+k)=0$,
into a form,
$\mathcal{A}_\mathcal{K}(E,Q;q) \, J_\mathcal{K}(N)
\equiv
\sum_{k\geq 0}   a_k(Q,q)  \, E^k
\, J_{\mathcal{K}}(N)=0
$,
where  the operators $Q$ and $E$ act on $J_\mathcal{K}(N)$ as
\begin{align*}
  E \, J_\mathcal{K}(N)
  & = J_\mathcal{K}(N+1)
  \\[2mm]
  Q \, J_\mathcal{K}(N)
  & = q^N \, J_\mathcal{K}(N)
\end{align*}
Then a claim of Ref.~\citen{SGarou03b} is that the A-polynomial
$A_\mathcal{K}(L,M)$ coincides with
\begin{equation}
  \label{AJ_identity}
  A_\mathcal{K}(L,M)
  =
  \mathcal{A}_\mathcal{K}(L,M^2; q=1)
\end{equation}
This AJ conjecture was checked with a help of \texttt{Mathematica}
package for trefoil, figure-eight knots~\cite{SGarou03b}, and
for  $5_2$,
$6_1$ knots~\cite{TTakata04a}.

For our case of the torus knot $\mathcal{T}_{s,t}$,
we can  easily check that
applying above procedure
the difference equations~\eqref{difference_Jones}
and~\eqref{torus_difference} reproduce the A-polynomial for the 
torus knots given in Ref.~\citen{XZhan03a}
(see also Refs.~\citen{CCGLS94a,PDShana00a});
\begin{equation*}
  A_{\mathcal{T}_{s,t}}(L, M)
  =
  \begin{cases}
    (L-1) \, (- 1 + L^2 \, M^{2 s t} ) 
    &
    \text{for $s,t>2$}
    \\[2mm]
    (L-1 ) \, ( 1 + L \,  M^{2(2m+1)} ) 
    &
    \text{for $(s,t)=(2,2m+1)$}
  \end{cases}
\end{equation*}
We can thus conclude that
\begin{theorem}
  \label{thm:AJ_torus}
  AJ conjecture~\eqref{AJ_identity} proposed in Ref.~\citen{SGarou03b}
  is true for the torus knots
  $\mathcal{K}=\mathcal{T}_{s,t}$.
\end{theorem}

\section{\mathversion{bold}
$q$-Hypergeometric Function and Colored Jones Polynomial}

At the end of this article, we comment on an explicit form of the
colored Jones polynomial for the torus knot.
Recalling a result of Ref.~\citen{KHikami02c},
we see that
the  colored Jones polynomial for the torus knot
$\mathcal{T}_{2,2m+1}$ can be written in a form of the
$q$-hypergeometric function.

Hereafter we use a standard notation of the $q$-product and the
$q$-binomial coefficient (see \emph{e.g.} Ref.~\citen{Andre76});
\begin{gather*}
  (x)_n = (x;q)_n = \prod_{k=1}^n (1 - x \, q^{k-1})
  \\[2mm]
  \begin{bmatrix}
    n \\
    m
  \end{bmatrix}_q
  =
  \frac{(q)_n}{
    (q)_{n-m} \, (q)_m
  }
\end{gather*}

\begin{theorem}[\cite{KHikami02c}]
  Let the function $H_{{2,2m+1}}(x)$ be defined by
  \begin{equation}
    \label{define_series}
    H_{{2,2m+1}}(x)
    =
    \sum_{
      k_m \geq \dots \geq k_2 \geq k_1 \geq 0
    }^\infty
    (x)_{k_m + 1} \, x^{k_m} \,
    \left(
      \prod_{i=1}^{m-1}
      q^{k_i (k_i+1)} \,
      x^{2 k_i} \,
      \begin{bmatrix}
        k_{i+1} \\
        k_i
      \end{bmatrix}_q
    \right)      
  \end{equation}
  for $|q|<1$ and $|x|<1$.
  Then we have
  \begin{equation}
    H_{{2,2m+1}}(x)
    =
    \sum_{n=0}^\infty \chi_{8m+4}(n) \,
    q^{\frac{n^2 - (2 m -1)^2}{8(2m+1)}} \,
    x^{\frac{n- (2m -1)}{2}}
  \end{equation}
  and it satisfies
  \begin{equation}
    \label{difference_H}
    H_{2,2m+1}(x)
    =
    1 - q \, x^2 -
    q^{2 m} \, x^{2 m +1} \,
    H_{2,2m+1}(q \, x)
  \end{equation}
\end{theorem}

For   $m=1$ case,
see also
Ref.~\citen{Andre76}~[Chap.~2, p.~29] and
Ref.~\citen{DZagie01a}.

With an expression~\eqref{define_series}, we can  take a limit $x\to q^{-N}$;
an infinite sum terminates into  a finite sum due to
$(q^{-N})_{k}=0$ for $k>N$.
Comparing eq.~\eqref{difference_K_2} with eq.~\eqref{difference_H}
we obtain the colored Jones polynomial for the torus knot
$\mathcal{T}_{2,2m+1}$
as follows;
\begin{prop}
  The $N$-colored Jones polynomial for the torus knot $\mathcal{T}_{2,2\,m+1}$ is given by
  \begin{equation}
    J_{\mathcal{T}_{2,2m+1}}(N)
    =
    q^{m(1-N)}
    \sum_{
      k_m \geq \dots \geq k_2 \geq k_1 \geq 0
    }^\infty
    (q^{1-N})_{k_m} \, q^{-N k_m} \,
    \left(
      \prod_{i=1}^{m-1}
      q^{k_i (k_i+1-2 N)} \,
      \begin{bmatrix}
        k_{i+1} \\
        k_i
      \end{bmatrix}_q
    \right)      
  \end{equation}
\end{prop}

We see that the colored Jones polynomial for the trefoil $m=1$
coincides with results in Refs.~\citen{KHabi02a,TQLe03a}.
It should be remarked that in Ref.~\citen{GMasb03a}
constructed was  the colored Jones
polynomial for the twist knot, which gives a different expression of the colored Jones
polynomial for the trefoil.
By applying a twisting formula given in
Ref.~\citen{GMasb03a}, we have another expression
of the colored Jones polynomial in the form of
\emph{cyclotomic expansion} in a sense of Ref.~\citen{KHabi02a} as follows;
\begin{prop}
  The $N$-colored Jones polynomial for the torus knot $\mathcal{T}_{2,2\,m+1}$ is written as
  \begin{multline}
    \label{another_torus}
    J_{\mathcal{T}_{2,2  m+1}}(N)
    =
    q^{ m \, (1-N^2)}
    \sum_{k_m \geq \dots \geq k_2 \geq k_1 \geq 0}^\infty
    \frac{
      (q^{1-N})_{k_m} \, (q^{1+N})_{k_m}
    }{
      (q)_{k_m}
    }
    \\
    \times
    \left(
      \prod_{i=1}^{m-1}
      q^{
        (k_i - k_m)
        (k_i - k_m-1)
      }\, 
      \begin{bmatrix}
        k_{i+1} \\
        k_i
      \end{bmatrix}_q
    \right) \,
  \end{multline}
\end{prop}

For the torus knot $\mathcal{T}_{s,t}$ with   $s, t >2$, we do
not know a
general  expression of the colored Jones polynomial in terms of the
$q$-hypergeometric function;
based on Refs.~\citen{KHikami03c,KHikami03b},
we only have
\begin{multline}
  J_{\mathcal{T}_{3,4}}(N)
  =
  q^{3 (1-N)} 
  \sum_{n=0}^\infty
  (q^{1-N})_n \, q^{-2 N n}
  \\
  \times
  \left(
    \sum_{k=0}^{\lfloor  (n-1)/2 \rfloor}
    q^{2 k (k+1-N)  + N}
    \begin{bmatrix}
      n \\
      2 \, k+1
    \end{bmatrix}_q
    +
    \sum_{k=0}^{\lfloor  n/2 \rfloor}
    q^{2 k (k+1-N) }
    \begin{bmatrix}
      n +1 \\
      2 \, k+1
    \end{bmatrix}_q
  \right)
\end{multline}




\section*{Acknowledgments}
The author would like to thank
S.~Garoufalidis,
J.~Kaneko,
A.~N.~Kirillov,
H.~Murakami,
T.~Takata,
Y.~Yokota,
D.~Zagier, and S.~Zwegers for discussions.
This work is supported in part by Grant-in-Aid for Young Scientists
from the Ministry of Education, Culture, Sports, Science and
Technology of Japan.


\end{document}